# NEGATIVE ASSOCIATION IN UNIFORM FORESTS AND CONNECTED GRAPHS

G. R. GRIMMETT AND S. N. WINKLER

ABSTRACT. We consider three probability measures on subsets of edges of a given finite graph $G$, namely those which govern, respectively, a uniform forest, a uniform spanning tree, and a uniform connected subgraph. A conjecture concerning the negative association of two edges is reviewed for a uniform forest, and a related conjecture is posed for a uniform connected subgraph. The former conjecture is verified numerically for all graphs $G$ having eight or fewer vertices, or having nine vertices and no more than eighteen edges, using a certain computer algorithm which is summarised in this paper. Negative association is known already to be valid for a uniform spanning tree. The three cases of uniform forest, uniform spanning tree, and uniform connected subgraph are special cases of a more general conjecture arising from the random-cluster model of statistical mechanics.

## 1. Three random subgraphs

Throughout this paper, $G = (V, E)$ denotes a finite labelled graph with vertex set $V$ and edge set $E$. An edge $e$ with endpoints $x, y$ is written $e = \langle x, y \rangle$. We assume that $G$ has neither loops nor multiple edges. We shall consider three probability measures on the set of subsets of $E$, and shall discuss certain results and conjectures concerning these measures. Since each such measure is a uniform measure on a given subset of $E$, each of our conclusions and conjectures may be expressed as a purely combinatorial statement.

The three measures are given as follows. Let $\mathcal{F}$ be the set of all subsets of $E$ which contain no cycle, noting that the elements of $\mathcal{F}$ are exactly the subgraphs of $G$ which are forests. Let $\mathcal{C}$ be the set of all subsets $C$ of $E$ such that $(V, C)$ is connected, and let $\mathcal{T} = \mathcal{F} \cap \mathcal{C}$ be the set of all spanning trees of $G$. We write $F$, $C$, $T$ (respectively) for elements of $\mathcal{F}$, $\mathcal{C}$, $\mathcal{T}$ chosen uniformly at random, and we call them a *uniform forest*, a *uniform connected subgraph*, and a *uniform spanning tree*, respectively. Note that all subgraphs of $G = (V, E)$ considered in this paper have the full vertex set $V$, and are thus said to be 'spanning'.

The uniform spanning tree has been studied extensively, see [**3**, **18**] for example. The uniform forests of this paper are different from those of [**3**] in that, in [**3**], the term *uniform spanning forest* denotes effectively the probability measure on infinite subsets of a given infinite graph $G$ which is obtained as the infinite-volume weak limit of the uniform spanning tree on a finite subgraph of $G$.

We remind the reader that, as summarised in Section 4, these three measures arise as weak limits of the random-cluster measure on $G$ in certain limiting regimes of the parameters of this measure.

The random subset $S$ of $E$ is said to be *edge-negatively-associated* if

$$\mathbb{P}(e \in S, \ f \in S) \leq \mathbb{P}(e \in S)\mathbb{P}(f \in S) \quad \text{for all } e, f \in E, \ e \neq f,$$

where $\mathbb{P}$ denotes the appropriate probability measure.

There is a more general notion of negative association, as follows. Random subsets of $E$ take values in the power set of $E$ which we denote by $\Omega_E = \{0,1\}^E$. An element $\omega = (\omega(e) : e \in E) \in \Omega_E$ is identified with the set $\eta(\omega) = \{e \in E : \omega(e) = 1\}$ of edges $e$ with $\omega(e) = 1$. An *event* is a subset of $\Omega_E$. For $E' \subseteq E$ and an event $A$, we say that $A$ is *defined on* $E'$ if: for all $\omega, \omega' \in \Omega_E$, if $\omega'(e) = \omega(e)$ for all $e \in E'$, then either $\omega, \omega' \in A$ or $\omega, \omega' \notin A$. The sample space $\Omega_E$ is a partially ordered set with order relation

$$\omega \leq \omega' \quad \text{if} \quad \omega(e) \leq \omega'(e) \text{ for all } e \in E.$$

An event $A$ is called *increasing* if $\omega' \in A$ whenever $\omega' \geq \omega$ and $\omega \in A$.

We call the random subset $S$ of $E$ *negatively associated* if

$$\mathbb{P}(S \in A \cap B) \leq \mathbb{P}(S \in A)\mathbb{P}(S \in B)$$

for all pairs $A$, $B$ of increasing events with the property that there exists $E' \subseteq E$ such that $A$ is defined in $E'$ and $B$ is defined on its complement $\overline{E'} = E \setminus E'$. For a general account of negative association and its inherent problems, see [20], which includes that part of the following conjecture concerning the uniform forest $F$.

**Conjecture 1.1.** *For all finite graphs $G = (V, E)$, the uniform forest $F$ and the uniform connected subgraph $C$ are edge-negatively-associated.*

A stronger version of this conjecture is that $F$ and $C$ are negatively associated in the general sense described above.

Since $F$ and $C$ are chosen uniformly at random, Conjecture 1.1 may be re-written in the form of two questions concerning subgraph counts. Let $V = \{1, 2, \ldots, n\}$, and let $K$ be the set of $N = \binom{n}{2}$ edges of the complete graph on the vertex set $V$. Let $\mathcal{G}(V)$ denote the set of all graphs on $V$ with neither loops nor multiple edges; since $\mathcal{G}(V)$ is in one–one correspondence with the set of all subsets of $K$, we shall take $\mathcal{G}(V)$ to be this power set. Thus we speak of subsets of $K$ as being graphs on $V$. Let $E \subseteq K$. For $X \subseteq E$, let $M^X = M^X(E)$ be the number of subsets $E'$ of $E$ with $E' \supseteq X$ such that the graph $(V, E')$ is connected. Edge-negative-association for connected subgraphs amounts to the inequality

(1.2) $$M^{\{e,f\}} M^{\varnothing} \leq M^e M^f \quad \text{for all } e, f \in E, \ e \neq f.$$

Here and later in this context, singleton sets are denoted without their braces, and any empty set is suppressed.

In the second such question, we ask if the same inequality is valid with $M^X$ redefined as the number of subsets $E'$ containing $X$ such that $(V, E')$ is a forest. (See [15, 20].)

With $E$ fixed as above, and with $X, Y \subseteq E$, let $M_Y^X = M_Y^X(E)$ denote the number of subsets $E'$ of $E$ of the required type such that $E' \supseteq X$ and $E' \cap Y = \varnothing$. Inequality (1.2) is easily seen to be equivalent to the inequality

(1.3) $$M^{\{e,f\}} M_{\{e,f\}} \leq M_f^e M_e^f.$$

The corresponding statement for the uniform spanning tree is known.

**Theorem 1.4.** *The uniform spanning tree $T$ is negatively associated.*

Thus was proved in [**7**], see also [**3**]. The lesser statement of edge-negative-association of $T$ is a consequence of the 1847 work of Kirchhoff [**16**] concerning electrical networks. See Section 3 for a brief account of the proof of edge-negative association for the uniform spanning tree.

The main advance reported in this note is the following.

**Theorem 1.5.** *If $G = (V, E)$ has eight or fewer vertices, or has nine vertices and eighteen or fewer edges, then the uniform forest $F$ has the edge-negative-association property.*

This we have verified by direct numerical calculation, using a combination of computer algorithms described in Section 5. That is, for every graph $G = (V, E)$ with eight or fewer vertices, and for $e, f \in E$, we have used a computer to count the numbers of forests appearing in (1.2), and we have thus checked that (1.2) is valid. This non-trivial calculation was achieved via a method summarised in Section 5 and expanded further in Sections 6 and 7. The central difficulty is that there is currently known no computationally efficient way of counting the number of forests contained in a given graph. Indeed, this problem is #P-complete and presumably computationally intractable, see [**1, 2, 21**] and also [**14**]. Thus, some cleverness is needed to complete the computations in a reasonable time. In contrast, counting the number of spanning trees of a graph is equivalent to the evaluation of a certain determinant, a computationally easy problem, see [**17**].

We point out that there are $2^{28}$ labelled graphs on eight vertices, and $2^{36}$ on nine vertices. Our computation for graphs of size eight required about 40 minutes on a PC running Linux, using a Pentium 3 processor with an official clock speed of 2.8GHz, and with approximately 2Gb RAM. Neither speed nor memory capacity was tested by the case $n = 8$. The case $n = 9$ proved much more demanding in computational terms, and exhausted the RAM after running for about 9 hours and completing the calculation for subgraphs having 16 or fewer edges. After taking certain measures which preserve RAM at the expense of speed, we were able to handle subgraphs with 18 edges or fewer in approximately 140 hours.

There are more general versions of Conjecture 1.1. First, one might ask whether the following holds instead of (1.2),

(1.6) $$M^{\{e,f\}}(\alpha) M^\varnothing(\alpha) \leq M^e(\alpha) M^f(\alpha) \quad \text{for } e \neq f,\ 0 < \alpha < \infty,$$

where

$$M^X(\alpha) = \sum_{\substack{E': X \subseteq E' \subseteq E \\ (V, E') \text{ has property } \Pi}} \alpha^{|E'|},$$

and $\Pi$ is either the property of either being a forest or the property of being connected. This reduces to (1.2) when $\alpha = 1$. We shall encounter such a generalisation in the context of the random-cluster measures in Section 4.

For a more general formulation, let $\mathbf{p} = (p_e : e \in K)$ be a collection of non-negative numbers indexed by $K$, and let $E \subseteq K$. For a subset $E' \subseteq E$, we write
$$f_{\mathbf{p}}(E') = \prod_{e \in E'} p_e.$$

We now ask whether (1.2) holds with $M^X = M^X(\mathbf{p})$ defined by

(1.7) $$M^X(\mathbf{p}) = \sum_{\substack{E': X \subseteq E' \subseteq E \\ (V, E') \text{ has property } \Pi}} f_{\mathbf{p}}(E'),$$

with $\Pi$ as above. Note that (1.2) becomes a polynomial inequality in $|E|$ real variables. We are grateful to Alan Sokal for his encouragement to formulate the problem in this manner. Such a formulation is natural when the problem is cast in the context of the Tutte polynomial, see Section 4.

We have discussed only graphs $G = (V, E)$ with neither loops nor multiple edges. No extra generality is gained in the consideration of edge-negative-association by introducing loops, since no forest contains loops, and in addition the inclusion or not of loops has no effect on whether a given subgraph is connected. The situation with multiple edges is more interesting, and is related to the general form of the conjecture of edge-negative-association discussed around (1.7).

## 2. Planar duality

The graph $G = (V, E)$ is called *planar* if it may be drawn in the plane in such a way that edges intersect only at their endpoints. Suppose that $G$ is planar and has been embedded in the plane in such a way. We write $G$ also for the planar embedding of $G$, and we obtain the dual graph $G^{\mathrm{d}} = (V^{\mathrm{d}}, E^{\mathrm{d}})$ as follows. We place a dual vertex within each face of $G$, including the infinite face. For each $e \in E$ we place a dual edge $e^{\mathrm{d}} = \langle x^{\mathrm{d}}, y^{\mathrm{d}} \rangle$ joining the two dual vertices lying in the two faces of $G$ abutting $e$; if these two faces are the same, then $x^{\mathrm{d}} = y^{\mathrm{d}}$ and $e^{\mathrm{d}}$ is a loop. Thus $E^{\mathrm{d}}$ is in one–one correspondence to $E$, and $G^{\mathrm{d}}$ is itself planar. Any configuration $\omega \in \Omega_E$ ($= \{0,1\}^E$) gives rise to a dual configuration $\omega^{\mathrm{d}}$ lying in the space $\Omega_E^{\mathrm{d}} = \{0,1\}^{E^{\mathrm{d}}}$ defined by $\omega^{\mathrm{d}}(e^{\mathrm{d}}) = 1 - \omega(e)$. As before, to each configuration $\omega^{\mathrm{d}}$ corresponds the set $\eta(\omega^{\mathrm{d}}) = \{e^{\mathrm{d}} \in E^{\mathrm{d}} : \omega^{\mathrm{d}}(e^{\mathrm{d}}) = 1\}$ of its 'open edges'. Note that a given planar graph may have more than planar embedding, and that different embeddings may have non-isomorphic dual graphs. It is easily seen that $G$ is the dual graph of $G^{\mathrm{d}}$.

**Theorem 2.1.** *A configuration $\omega \in \Omega$ corresponds to a forest of $G$ if and only if the dual configuration $\omega^{\mathrm{d}}$ corresponds to a connected subgraph of $G^{\mathrm{d}}$.*

*Proof.* If $\eta(\omega)$ contains a cycle $\gamma$, then the dual vertices corresponding to faces of $G$ inside $\gamma$ are disconnected from those corresponding to faces outside $\gamma$. Thus the graph $(V^{\mathrm{d}}, \eta(\omega^{\mathrm{d}}))$ is disconnected. It is easy also to see the converse: if $\eta(\omega)$ contains no cycle, then $(V^{\mathrm{d}}, \eta(\omega^{\mathrm{d}}))$ is connected. $\square$

Let $X, Y \subseteq E$. We deduce from Theorem 2.1 and its proof that the number of forests of the planar graph $G$ containing every edge in $X$ but no edge of $Y$ is equal to the number of connected subgraphs of $G^{\mathrm{d}}$ containing every edge of $Y^{\mathrm{d}} = \{f^{\mathrm{d}} : f \in Y\}$ but no edge in $X^{\mathrm{d}} = \{f^{\mathrm{d}} : f \in X\}$. It follows by the equivalence of (1.2) and (1.3) that, for any planar graph $G$, the uniform spanning forest $F$ is edge-negatively-associated if and only if the uniform connected subgraph $C^{\mathrm{d}}$ of $G^{\mathrm{d}}$ has this property. In particular, Theorem 1.5 implies the edge-negative-association of the uniform connected subgraph of any planar graph with eight or fewer faces (including the infinite face) whose dual has no multiple edges.

## 3. Edge-negative-association of uniform spanning tree

Since the emphasis of this note is upon the property of edge-negative-association, we include a very short account of its proof for the uniform spanning tree. Kirchhoff [16] showed how to relate the current which flows along the edge $e$ in an electrical network to certain counts of spanning trees. Consider an electrical network on the connected graph $G$ in which each edge corresponds to a unit resistor. The relevant fact from the theory of electrical networks is that, if a unit current flows from a *source* vertex $s$ to a *sink* vertex $t$, then the current which flows along the edge $e = \langle x, y \rangle$ in the direction $xy$ equals $N(s, x, y, t)/N$, where $N$ is the number of spanning trees of $G$ and $N(s, u, v, t)$ is the number of spanning trees whose unique path from $s$ to $t$ passes along the edge $\langle u, v \rangle$ in the direction $uv$.

Let $e = \langle x, y \rangle$. By the above, $\mathbb{P}(e \in T)$ equals the current flowing along $e$ when a unit current flows through $G$ from source $x$ to sink $y$. By Ohm's Law, this equals the potential difference between $x$ and $y$, which in turn equals the effective resistance $R_G(x, y)$ of the network between $x$ and $y$.

Let $f \in E$, $f \neq e$, and denote by $G.f$ the graph obtained from $G$ by contracting the edge $f$. There is a one–one correspondence between spanning trees of $G.f$ and spanning trees of $G$ containing $f$. Therefore, $\mathbb{P}(e \in T \mid f \in T)$ equals the effective resistance $R_{G.f}(x, y)$ of the network $G.f$ between $x$ and $y$. The so-called Rayleigh principle states that effective resistance is a non-decreasing function of the individual edge-resistances, and it follows as required that $R_{G.f}(x, y) \leq R_G(x, y)$.

We note that the usual proof of the Rayleigh principle makes use of the Thomson/Dirichlet variational principle, which in turn asserts that, amongst all unit flows from source to sink, the true current flow of size one is that which minimises the dissipated energy.

A good account of the Kirchhoff theorem, above, may be found in [4]. Further accounts of the mathematics of electrical networks include [6] and [18], the latter containing also much material about the uniform spanning tree.

## 4. The random-cluster measure

Let $\Omega_E = \{0, 1\}^E$ be the set of 0/1-vectors indexed by $E$, and let $\omega \in \Omega_E$. An edge $e$ is termed *open* in $\omega$ if $\omega(e) = 1$, and is termed *closed* otherwise. We write $\eta(\omega) = \{e \in E : \omega(e) = 1\}$ for the set of open edges of $\omega$, and $k(\omega)$ for the number of components of the open graph $(V, \eta(\omega))$. Amongst the various probability measures on $\Omega_E$, probably the most studied is product measure $\mu_p$ with some given density $p$. The subsequent theory is usually referred to either as 'percolation' when the graph is part of a lattice, or as 'random graphs' when $G$ is complete; see [5, 9, 13]. An important extension of product measure is the so-called random-cluster measure given as follows. Let $0 \le p \le 1$ and $q > 0$. We define the probability measure $\phi_{p,q}$ on $\Omega_E$ by

$$\phi_{p,q}(\omega) = \frac{1}{Z(p,q)} \left\{ \prod_{e \in E} p^{\omega(e)}(1-p)^{1-\omega(e)} \right\} q^{k(\omega)}, \qquad \omega \in \Omega_E,$$

where the 'partition function' $Z(p,q)$ is given by

$$Z(p,q) = \sum_{\omega \in \Omega} \left\{ \prod_{e \in E} p^{\omega(e)}(1-p)^{1-\omega(e)} \right\} q^{k(\omega)}.$$

This measure was introduced around 1970 by Fortuin and Kasteleyn, and has proved very useful in the study of Ising and Potts models. When $q = 1$, $\phi_{p,q}$ is simply product measure with density $p$. When $q \in \{2, 3, \dots\}$, $\phi_{p,q}$ corresponds in a certain way to a Potts model on $G$ with $q$ states available at each vertex. It is a very useful property that $\phi_{p,q}$ is positively associated when $q \ge 1$, and this property does not hold when $q < 1$ (except for graphs with no cycles). Only little is known about the measure $\phi_{p,q}$ when $q < 1$, although a few facts can be proved. See [10] for a history and a review of the measure, and [11] for a more extensive account with proofs.

We show next as in [10] that the uniform forest, uniform connected graph and uniform spanning tree are obtained as weak limits of $\phi_{p,q}$ as $p, q \downarrow 0$. This was known to Fortuin and Kasteleyn for spanning trees, see [8, 12], and to Pemantle [20] and perhaps others for forests.

Suppose that $0 < p < 1$, and consider the weak limit of $\phi_{p,q}$ as $q \downarrow 0$. Because of the term $q^{k(\omega)}$ in the definition of $\phi_{p,q}$, all the mass is concentrated in the limit on those configurations $\omega$ for which $k(\omega)$ is a minimum, that is, with $k(\omega) = 1$. Since the 'weight' of a configuration $\omega$ is proportional to $\{p/(1-p)\}^{|\eta(\omega)|} q^{k(\omega)}$, the limiting probability measure is given, with $\beta = p/(1-p)$, by

$$\mu_{\mathcal{C},\beta}(\omega) = \frac{1}{Z_{\mathcal{C},\beta}} \beta^{|\eta(\omega)|}, \qquad \omega \in \Omega_{\mathcal{C}},$$

where $\Omega_{\mathcal{C}}$ is the set of all $\omega$ such that $(V, \eta(\omega)) \in \mathcal{C}$, the set of connected subgraphs, and

$$Z_{\mathcal{C},\beta} = \sum_{\omega \in \Omega_{\mathcal{C}}} \beta^{|\eta(\omega)|}.$$

Thus $\mu_{\mathcal{C},\beta}$ is concentrated on the set of connected subgraphs of $G$, and $\mu_{\mathcal{C},1}$ is uniform measure on $\mathcal{C}$.

Suppose secondly that $p = p_q$ is related to $q$ in such a way that $p \to 0$ and $q/p \to 0$ as $q \to 0$. The partition function $Z(p,q)$ may be expressed in the form

$$Z(p,q) = (1-p)^{|E|} \sum_{\omega \in \Omega} \left( \frac{p}{1-p} \right)^{|\eta(\omega)|+k(\omega)} \left( \frac{q(1-p)}{p} \right)^{k(\omega)}.$$

Note that $p/(1-p) \to 0$ and $q(1-p)/p \to 0$ as $q \to 0$. Now $k(\omega) \geq 1$ and $|\eta(\omega)| + k(\omega) \geq |V|$ for all $\omega \in \Omega$; these two inequalities are satisfied simultaneously with equality if and only if $\eta(\omega)$ is a spanning tree of $G$. It follows that, in the limit as $q \to 0$, the 'mass' is concentrated on such configurations, and it is easily seen that the limit mass is uniformly distributed. That is, $\lim_{q \downarrow 0} \phi_{p,q}$ is a probability measure which selects, uniformly at random, a spanning tree of $G$.

Next we suppose $p = q$ and consider the limit as $q \downarrow 0$. The important component in the weight of $\omega$ is $q^{|\eta(\omega)|+k(\omega)}$. Now $|\eta(\omega)| + k(\omega) \geq |V|$, with equality if and only if the set of open edges contains no cycle. It follows that $\phi_{p,q}$ converges weakly as $q \downarrow 0$ to the probability measure

$$\mu_{\mathcal{F},1}(\omega) = \frac{1}{|\Omega_{\mathcal{F}}|}, \qquad \omega \in \Omega_{\mathcal{F}},$$

where $\Omega_{\mathcal{F}}$ is the set of all $\omega$ such that $\eta(\omega) \in \mathcal{F}$, the set of forests in $G$. Thus $\mu_{\mathcal{F},1}$ is uniform measure on $\mathcal{F}$. More generally, take $p = \alpha q$ where $\alpha \in (0, \infty)$ is constant, and take the limit as $q \downarrow 0$. The limiting measure is now

$$\mu_{\mathcal{F},\alpha}(\omega) = \frac{1}{Z_{\mathcal{F},\alpha}} \alpha^{|\eta(\omega)|}, \qquad \omega \in \Omega_{\mathcal{F}},$$

where

$$Z_{\mathcal{F},\alpha} = \sum_{\omega \in \Omega_{\mathcal{F}}} \alpha^{|\eta(\omega)|}.$$

It may be conjectured that $\phi_{p,q}$ is negatively associated in the sense of [7] when $q < 1$, and in particular that $\phi_{p,q}$ is edge-negatively-associated in that

$$\phi_{p,q}(e \text{ is open}, f \text{ is open}) \leq \phi_{p,q}(e \text{ is open})\phi_{p,q}(f \text{ is open}),$$

for distinct edges $e$, $f$. A verification of this conjecture would be valuable in studying the random-cluster measure with $q < 1$.

The partition function is an evaluation of the Tutte polynomial of the graph $G$, and for this reason and others there is a strong relationship between the random-cluster model and such polynomials. We do not describe this here, but refer the reader to the references in [10, 11].

## 5. The algorithm for counting forests

We summarise next our approach to the computations necessary for Theorem 1.3. The computer code and some additional documentation is available at the web page `http://www.statslab.cam.ac.uk/~grg/usf/`.

To verify (1.3) for all edge sets $E$ on $n$ vertices is a numerical problem of considerable complexity. We give a brief summary of such an algorithm which, by dint of careful book-keeping, confirms (1.3) whenever $n \leq 8$. Our general approach is sketched in this section, and elaborated in Sections 6 and 7.

5.1 OVERVIEW

Let $V = \{1, 2, \ldots, n\}$ where $n \geq 3$, and let $K$ be the set of $N = \binom{n}{2}$ edges of the complete graph on $V$, as before. The edge $e$ of $K$ with endpoints $x$ and $y$ such that $x < y$ is written $e = \langle x, y \rangle$. There is a total (lexicographic) order on $K$ given by $\langle u, v \rangle < \langle x, y \rangle$ if either $u < x$ or $u = x$, $v < y$. For $X, Y \subseteq E$, we write $\mathcal{F}_Y^X(E)$ for the set of all forests $F \ (\subseteq E)$ such that $X \subseteq F$ and $F \cap Y = \varnothing$, and we refer to $X$ and $Y$ as the 'constraints' of such forests. We write $M_Y^X(E) = |\mathcal{F}_Y^X(E)|$. If either $X$ or $Y$ is empty, we omit it from the notation. We shall sometimes write $A + B$ for the union of two sets $A$, $B$, and $A - B$ for $A \setminus B$, and we shall generally omit braces around singleton sets when no ambiguity is thus introduced.

It is required to prove (1.3), namely,

$$(5.1) \qquad M^{\{e,f\}}(E) M_{\{e,f\}}(E) \leq M_f^e(E) M_e^f(E)$$

for all $e, f \in E$, $e \neq f$, and for all $E \subseteq K$. That is, we seek to validate (5.1) for all suitable triples $(E, e, f)$, and we refer to such a triple as a 'conjecture instance'. Note that conjecture instances may be classified into two categories, depending on whether or not $e$ and $f$ share a vertex. With this in mind, we set $E_1 = \{e_1, f_1\} = \{\langle 1, 2 \rangle, \langle 1, 3 \rangle\}$ and $E_2 = \{e_2, f_2\} = \{\langle 1, 2 \rangle, \langle 3, 4 \rangle\}$, each such set being representative of one of the two categories.

The general approach is as follows. For each $s = 2, 3, \ldots, N$ in order, we generate the family $\mathcal{E}_s$ of all subsets of $K$ with cardinality $s$ and containing either $E_1$ or $E_2$. For each $s$, $i$, and each $E \in \mathcal{E}_s$ satisfying $E_i \subseteq E$, we check (5.1) for the conjecture instance $(E, e_i, f_i)$. To do so requires computations of the four numbers $M_Y^{E_i - Y}(E)$ for $Y = \varnothing, \{e_i\}, \{f_i\}, E_i$. We expand on such computations in Subsection 5.2.

This method, if implemented naively, would be computationally infeasible in practice, since the size of the problem grows very quickly with $n$, and since no efficient method is known for calculating quantities of the form $M_Y^X(E)$. We propose two substantial economies of scale. The first aims to reduce the number of conjecture instances to be studied, (a) by eliminating instances which are 'isomorphic' to an instance already resolved, and (b) by developing a criterion for deciding when two conjecture instances are effectively equivalent, in the sense that (5.1) holds for the second whenever it holds for the first. This is discussed in much more detail in Section 7.

Our second economy pertains to the calculation of the terms $M_Y^X(E)$. As described in Subsection 5.4, such a calculation involves many sub-calculations, each of which may have been encountered earlier in the same or some related form. We aim to keep a record of certain earlier calculations, the better to avoid duplication of computational effort.

Each of these measures hinges on the definition of a suitable 'index function', by use of which one may recognise when a proposed calculation is unnecessary; this is sketched in Subsection 5.3. Our index function is presented in Section 6, and is based on the `nauty` package developed by Brendan McKay in the context of the graph-isomorphism problem.

## 5.2 Evaluating $M_Y^X(E)$

We refer to the triple $(E, X, Y)$ as a 'counting problem', and we shall need to compute $M_Y^X(E)$ for general counting problems. Note first that

$$(5.2) \qquad M_Y^X(E) = M^X(E - Y),$$

whence it suffices to compute certain $M^X(D)$. A brute-force strategy for calculating $M_Y^X(E)$ is to enumerate all valid forests $F$ of the same size before proceeding to the next layer in the hierarchy of graphs sandwiched between $X$ and $E - Y$. Let $\mathcal{F}_s = \{F \in \mathcal{F}^X(E - Y) : |F| = s\}$ denote the layer of valid forests of size $s$ for $s = |X|, \ldots, |E - Y|$. The $\mathcal{F}_s$ partition the set $\mathcal{F}^X(E - Y)$.

The construction of the layers $\mathcal{F}_s$ is described inductively, noting that $\mathcal{F}_{|X|} = \{X\}$. Having found $\mathcal{F}_s$, we construct $\mathcal{F}_{s+1}$ by considering every $F \in \mathcal{F}_s$ in turn. The forest $F$ can be thought of as a (vertex- and edge-)disjoint union of connected components $C_i$. We seek to enlarge $F$ by one edge without forming a cycle. Every such augmentation arises by merging two distinct components (whereas adding an edge to any given component must form a cycle). Thus, from every pair $C_i$, $C_j$ of distinct components, we obtain new elements $F'$ of $\mathcal{F}_{s+1}$ by adding an edge joining $C_i$ to $C_j$.

We have now that

$$(5.3) \qquad M_Y^X(E) = M^X(E - Y) = \sum_{s=|X|}^{|E-Y|} |\mathcal{F}_s|.$$

## 5.3 Indexing the counting problems

The ensuing algorithm is very time-consuming when $X$ is small and $E - Y$ is large, since the layers $\mathcal{F}_s$ grow fast in $n$ and $s$. It is therefore necessary to reduce to a minimum the number of its calls. We do this by devising some 'index function' $i$ that indexes classes of counting problems $(E, X, Y)$ having a common value of $M_Y^X(E)$. Let this value be $M_\iota$ for all problems with index satisfying $i(E, X, Y) = \iota$. Thus, for each $\iota$, one computation of $M_\iota$ suffices; the obtained value may be buffered in a database and retrieved to give $M_Y^X$ for all problems $(E, X, Y)$ with the same index $\iota$.

There are many possible choices of such a function $i$, and we seek one whose range is small and whose values are easy to compute. We show in the next section how such a function may be defined with the following key properties. Let $(E, X, Y)$ and $(E', X', Y')$ be two counting problems such that $i(E, X, Y) = i(E', X', Y')$. Then, in ways to be made more precise in Section 6,

(K1) $X$ and $X'$ as well as $E - (X + Y)$ and $E' - (X' + Y')$ have similar structure, respectively, and

(K2) $X$ and $E - (X + Y)$ are composed in a way similar to $X'$ and $E' - (X' + Y')$.

Properties (K1)–(K2) imply as required that $M^X(E - Y) = M^{X'}(E' - Y')$.

## 5.4 Recursive formula for $M_Y^X(D)$

A large number of explicit computations of $M^X(D)$ for $X \subseteq D \subseteq E$ remain inevitable. The cost of computing $M^X(D)$ can be large when $X$ is small, and grows very fast in the size of $D$. It is therefore preferable to expand this quantity into a sum over certain $M^{X'}(D')$ for larger sets $X'$ of constraints and smaller edge sets $D'$. These may be faster to compute, or may even be available from earlier computations. By partitioning the set of forests of $D$ according to whether or not they contain a given edge $d \in D - X$, we see that

$$(5.4) \qquad M^X(D) = M^X(D - d) + M^{X+d}(D) \qquad \text{for } d \in D \setminus X.$$

By assumption, $X \subseteq D$, and we write $D - X = \{d_1, d_2, \ldots, d_k\}$. Let $X_j = X + \{d_1, d_2, \ldots, d_{j-1}\}$ and $D_j = D - d_j$. By iteration of (5.4),

$$(5.5) \qquad M^X(D) = \sum_{j=1}^k M^{X_j}(D_j) + M^D(D),$$

a formula which substantially aids numerical economy.

## 6. Canonical labelling

We sketch next how to label counting problems in order to build up a database of problems considered so far. We make substantial use of concepts and methods introduced by Brendan McKay in [**19**] and implemented in his `nauty` package to be found at http://cs.anu.edu.au/~bdm/nauty/, and we begin with some terms.

### 6.1 Partitions

A *partition* $\pi = (V_1, V_2, \ldots, V_k)$ of the set $V$ is a sequence of disjoint non-empty subsets of $V$ whose union is $V$. The set of all partitions of $V$ will be denoted by $\Pi(V)$. If $V$ is the vertex set of a graph $G = (V, E)$, and $\pi \in \Pi(V)$, then the pair $(G, \pi)$ is termed a *partitioned graph*. The components $V_i$ of a partition $\pi$ are called its *cells*. Note that the order of the cells is significant, but the order of the vertices within each cell is not. A *singleton cell* is a cell $V_i = \{v\}$ with a single element $v$, which is said to be *fixed* by $\pi$. If every cell of $\pi$ is trivial, then $\pi$ is called *discrete*. If $\pi$ contains only one cell $V_1 = V$, then $\pi$ is called the *unit* partition. A partition $\pi$ with $k$ cells may be identified with a vertex-colouring $c : V \to \{1, 2, \ldots, k\}$ by assigning colour $i$ to all vertices in the $i$th cell $V_i$. In this alternative language, cells are also called *colour classes*.

### 6.2 Permutations

Let $\gamma$ be a permutation of $V$. The image of $v \in V$ under $\gamma$ will be denoted $v^\gamma$. Similarly, for $U \subseteq V$, we write $U^\gamma = \{u^\gamma : u \in U\}$. If $\pi = (V_1, V_2, \ldots, V_k) \in \Pi(V)$, we set $\pi^\gamma = (V_1^\gamma, V_2^\gamma, \ldots, V_k^\gamma)$. If $G = (V, E)$ is a graph, then $G^\gamma$ denotes the graph in which vertices $v^\gamma$ and $w^\gamma$ are adjacent if and only if $v$ and $w$ are adjacent in $G$. Since a permutation $\gamma$ acts on unordered pairs $\{u, v\}$ of vertices, it acts also on edges. We now see that $G^\gamma = (V, E^\gamma)$.

Two partitions $\pi_1, \pi_2 \in \Pi(V)$ are called *compatible* if there exists a permutation $\gamma$ of $V$ such that $\pi_2 = \pi_1^\gamma$. Thus $\pi_1$ and $\pi_2$ are compatible if and only if the vectors of cardinalities of their cells are the same. Two graphs $G_1, G_2 \in \mathcal{G}(V)$ are said to be *isomorphic*, written $G_1 \sim G_2$, if there exists a permutation $\gamma$ of $V$ such that $G_2 = G_1^\gamma$. If this holds, we call $\gamma$ an *isomorphism* of $G_1$ to $G_2$. If the graphs $G_i$ carry partitions $\pi_i$, we call the isomorphism $\gamma$ *partition-preserving* if, in addition, $\pi_2 = \pi_1^\gamma$.

The *automorphism group* $\mathrm{Aut}(G, \pi)$ of a partitioned graph $(G, \pi)$ is the set of all permutations $\gamma$ such that $G^\gamma = G$ and $\pi^\gamma = \pi$. Since the order of cells in partitions is significant, the latter condition requires that $\gamma$ fixes the cells in $\pi$ setwise. If $\pi$ is the unit partition, $\mathrm{Aut}(G, \pi)$ is the usual automorphism group $\mathrm{Aut}(G)$ of $G$.

### 6.3 Canonical labelling

A *canonical labelling map* is a function $\mathcal{C} : \mathcal{G}(V) \times \Pi(V) \to \mathcal{G}(V)$ such that, for any graph $G = (V, E)$, any partition $\pi$ of $V$, and any permutation $\gamma$ of $V$, we have
(C1) $\mathcal{C}(G, \pi) \sim G$,
(C2) $\mathcal{C}(G^\gamma, \pi^\gamma) = \mathcal{C}(G, \pi)$, and
(C3) if $\mathcal{C}(G, \pi^\gamma) = \mathcal{C}(G, \pi)$, there exists $\delta \in \mathrm{Aut}(G)$ such that $\pi^\gamma = \pi^\delta$.
It was established in [19] that a canonical labelling map exists, and such a map may be found which is generally easy to compute. Henceforth, we assume that $\mathcal{C}$ is such a map. If $\pi = (V)$ is the unit partition, we may write $\mathcal{C}(G)$ for $\mathcal{C}(G, (V))$. For $G = (V, E)$ and a partition $\pi$, we sometimes write $\mathcal{C}(E, \pi)$ for $\mathcal{C}(G, \pi)$.

The main use of a canonical labelling is to solve certain graph isomorphism problems, as illustrated by the following theorem which we shall use later.

**Theorem 6.1 [19, Thm 2.2].** *Let $G_1, G_2 \in \mathcal{G}(V)$ be graphs, $\pi_1, \pi_2 \in \Pi(V)$ compatible partitions, $\gamma$ a permutation of $V$, and let $\mathcal{C}$ be a canonical labelling map. Then $\mathcal{C}(G_1, \pi_1) = \mathcal{C}(G_2, \pi_2)$ if and only if there exists a partition-preserving isomorphism from $G_1$ to $G_2$.*

### 6.4 Sufficient condition for $M^X(D) = M^{X'}(D')$

We shall next define the index function $i$ referred to in Section 5. Let $D, D' \subseteq K$ and let $X \subseteq D$, $X' \subseteq D'$. We require a rather rigid criterion on the composition of the pairs $X, D - X$ and $X', D' - X'$ in order to guarantee that $M^X(D) = M^{X'}(D')$. The condition that we shall present amounts to requiring that the pairs $X, X'$ and $D - X, D' - X'$ are each isomorphic, and in addition that $X, D - X$ and $X', D' - X'$ 'fit together' in the same manner. It will be convenient to achieve this by use of the notion of a partition, and of McKay's canonical labelling map $\mathcal{C}$.

Let $\mathcal{C}(D)$ denote the canonically-labelled isomorph of $D$ (with unit partition), so that $\mathcal{C}(D) = D^\delta$ for some permutation $\delta$ of $V$; we fix such $\delta = \delta(D)$ henceforth. Given two edge sets $A, B \subseteq K$, we define their *interface* $\mathcal{I}(A, B)$ to be the set of vertices $v$ such that there exist $a \in A$, $b \in B$ with both $a$ and $b$ incident with $v$. Elements of $\mathcal{I}(A, B)$ are called *interface vertices*. Thus $\mathcal{I}(X, D - X)$ represents the interface of the complementary subgraphs $X$ and $D - X$ of $D$.

We wish to 'fix' this interface whilst allowing portions supported by vertices 'inside' $X$ and $D - X$ to vary up to isomorphism. For this purpose, we colour each vertex $v \in \mathcal{I}(X, D - X)$ with its label in the canonical isomorph of $D$, namely $v^\delta$. This induces a colouring $c : V \to \{0, 1, \ldots, n\}$ of the vertex set of $D$, each vertex $v$ taking colour $v^\delta$ if $v \in \mathcal{I}(X, D - X)$ and $0$ (black) otherwise.

The colouring $c$ gives rise to a partition of $V$ as follows. For $0 \leq i \leq n$, we set $V_i = \{v \in V : c(v) = i\}$. If $i \neq 0$, then $V_i$ is either empty or a singleton cell. Dropping all empty cells and relabelling the others while retaining the original order, we obtain a partition denoted $\pi = \pi(X, D - X, \delta)$ of the form $(V_0, \{c_1\}, \ldots, \{c_k\})$, where $k = |\mathcal{I}(X, D - X)|$.

The required index function $i$ is given in the next theorem.

**Theorem 6.2.** *Let $(E, X, Y)$ be a counting problem, let $D = E - Y$ and $\pi = \pi(X, D - X, \delta)$ where $\delta = \delta(D)$. Define the index function*

$$i(E, X, Y) = \big(|\mathcal{I}(X, D - X)|, \mathcal{C}(X, \pi), \mathcal{C}(D - X, \pi)\big).$$

*If $i(E, X, Y) = i(E', X', Y')$ then $M_Y^X(E) = M_{Y'}^{X'}(E')$.*

Since $M_Y^X(E) = M^X(E - Y)$, see (5.2), we shall have recourse to indices of the form $i(E - Y, X, \varnothing)$ which we write as $j(E - Y, X)$. Thus

(6.3) $$j(E', X) = i(E', X, \varnothing).$$

*Proof.* Assume that $i(E, X, Y) = i(E', X', Y')$. Since $|\mathcal{I}(X, D - X)| = |\mathcal{I}(X', D' - X')|$, the partitions $\pi$ and $\pi'$ are compatible, and thus we may apply Theorem 6.1 to the pairs $(X, \pi)$, $(X', \pi')$, and $(D - X, \pi)$, $(D' - X', \pi')$. We deduce from the equality of $i(E, X, Y)$ and $i(E', X', Y')$, and Theorem 6.1, that there exist partition-preserving permutations $\rho : (X, \pi) \mapsto (X', \pi')$ and $\sigma : (D - X, \pi) \mapsto (D' - X', \pi')$. Thus

(6.4) $$\begin{aligned} X' = X^\rho, \quad \pi' = \pi^\rho, \\ D' - X' = (D - X)^\sigma, \quad \pi' = \pi^\sigma. \end{aligned}$$

Since $\pi' = \pi^\rho = \pi^\sigma$, and since each point in $\mathcal{I}(X', D' - X')$ is a singleton cell of $\pi'$, we have that the permutations $\rho$ and $\sigma$ agree on the interface $\mathcal{I}(X, D - X)$.

For $F \subseteq K$, we write $V(F)$ for the set of endpoints of edges in $F$. Thus, $\overline{V(F)} = V - V(F)$ is the set of isolated vertices of the graph $(V, F)$. We define $\eta : V \to V$ by

(6.5) $$v^\eta = \begin{cases} v^\rho = v^\sigma & \text{if } v \in \mathcal{I}(X, D - X), \\ v^\rho & \text{if } v \in V(X) - \mathcal{I}(X, D - X), \\ v^\sigma & \text{otherwise,} \end{cases}$$

and we claim that $\eta$ is an isomorphism from $(V, D)$ to $(V, D')$. We check three facts in order to verify this.

(i) We claim that $\eta$ is one–one. Suppose on the converse that $x \neq y$ and $x^\eta = y^\eta$. Since $\rho$ and $\sigma$ are permutations, we may label $x$, $y$ in such a way that $x \in V(X) - \mathcal{I}(X, D - X)$ and $y \in V - V(X)$. Since $x \in V(X)$, we have that $x^\eta \in V(X')$. Since $y \in V - V(X)$, we have that $y^\eta \in V(D' - X') + \overline{V(D')}$. Since $x^\eta = y^\eta$, we deduce that $x^\eta = y^\eta \in \mathcal{I}(X', D' - X')$. Now $\rho$, $\sigma$ are partition-preserving, whence $x = y \in \mathcal{I}(X, D - X)$, a contradiction.

(ii) We claim that $\eta$ maps edges to edges. Let $\langle x, y \rangle \in D$. Either $\langle x, y \rangle \in X$, in which case $\langle x, y \rangle^\eta = \langle x^\rho, y^\rho \rangle \in X'$, or $\langle x, y \rangle \in D - X$, in which case $\langle x, y \rangle^\eta = \langle x^\sigma, y^\sigma \rangle \in D' - X'$. In either case, $\langle x^\eta, y^\eta \rangle \in D'$.

(iii) We claim that $\eta$ maps non-edges to non-edges. Assume that $\langle x^\eta, y^\eta \rangle \in D'$. If $\langle x^\eta, y^\eta \rangle \in X'$, then $\langle x^\eta, y^\eta \rangle = \langle x^\rho, y^\rho \rangle$ whence $\langle x, y \rangle \in X$. A similar argument holds with $X$, $\rho$ replaced by $D - X$, $\sigma$.

In summary, $\eta$ is a permutation of $V$ such that $e \in D$ if and only is $e^\eta \in D'$. Therefore, $\eta$ is an isomorphism.

It is easily deduced that $\eta$ induces a one–one correspondence between forests of $D$ containing $X$ and forests of $D'$ containing $X'$, whence $M^X(D) = M^{X'}(D')$ as required. $\square$

We make a peripheral remark regarding efficiency. In the above construction of the partition $\pi = \pi(X, D - X, \delta)$, it is not relevant for the validity of Theorem 6.2 to colour each vertex $v \in \mathcal{I}(X, D - X)$ precisely with its label $v^\delta$ in the canonical isomorph $\delta = \delta(D)$ of $D$. In principle, we could have chosen arbitrary distinct colours, thus removing the dependence of $\pi$ on $\delta$.

However, these colours determine the order in which interface vertices appear as singleton cells in the induced partition $\pi$. This order is in turn significant when applying the canonical labelling map $\mathcal{C}$. By consistently choosing the colour $v^\delta$ we maintain the same order across the whole class of isomorphs of $D$. As a consequence, $i$ sub-divides the space of counting problems $(E, X, Y)$ into larger classes, each with a common value $i(E, X, Y)$.

One drawback of the fact that the partition $\pi$ depends on $\delta(D)$ is that we have to recompute each of the three components of $i$ whenever a new graph $D$ is considered and the corresponding $\delta(D)$ affects $\pi$, even when $X$ or $D - X$ have not themselves changed. Whether or not the net effect is an improved performance will therefore depend essentially on the speed of evaluation of the canonical labelling map $\mathcal{C}$.

## 7. Redundancy in conjecture instances

We are required to check (5.1) for all conjecture instances $(E, e, f)$. Many such instances are isomorphic to one another, and thus it becomes necessary to devise a criterion for deciding, given two instances one of which has already been resolved, whether (5.1) holds automatically for the second.

We call the conjecture instances $(E, e, f)$ and $(E', e', f')$ *equivalent*, and write $(E, e, f) \sim (E', e', f')$, if the left and right hand sides of (5.1) are invariant under interchange of both sets of parameters, that is, if

$$
\begin{aligned}
M^{\{e,f\}}(E) M_{\{e,f\}}(E) &= M^{\{e',f'\}}(E') M_{\{e',f'\}}(E'), \quad \text{and} \\
M_f^e(E) M_e^f(E) &= M_{f'}^{e'}(E') M_{e'}^{f'}(E').
\end{aligned}
\tag{7.1}
$$

We propose two measures to reduce complexity: (a) to make a selection of conjecture instances which *ab initio* avoids many equivalent instances, and (b) to skip instances identified as equivalent to earlier ones, using the index function $i$ of Section 6.

7.1 EFFICIENT SELECTION

We begin by fixing $e$ and $f$ in one of the two categories of conjecture instances. Thus we take $E_1 = \{e_1, f_1\} = \{\langle 1,2\rangle, \langle 1,3\rangle\}$ and $E_2 = \{e_2, f_2\} = \{\langle 1,2\rangle, \langle 3,4\rangle\}$. It clearly suffices to verify (5.1) for all $E$ containing $E_1$, respectively $E_2$.

Let $\mathcal{E}_s = \{E \subseteq K : |E| = s,\ E_i \subseteq E \text{ for some } i\}$, $2 \leq s \leq N$. We shall not work with *every* member of $\mathcal{E}_s$ for each $s$, since the size of $\mathcal{E}_s$ grows exponentially with $n$ and $s$, and in addition many of its members give rise to equivalent conjecture instances. Instead we shall work with a subset $\mathcal{E}'_s \subseteq \mathcal{E}_s$ defined as follows.

For $E \subseteq K$, let $v_{\min} = v_{\min}(E)$ be the vertex having the least label which is isolated in $E$. If no vertex is isolated, we set $v_{\min} = n$. We say that $E$ has the property $\Pi_{\min}$ if each of the vertices $v_{\min}, v_{\min}+1, \ldots, n-1, n$ is isolated in $E$.

We set $\mathcal{E}'_2 = \mathcal{E}_2 = \{E_1, E_2\}$, and we make $\mathcal{E}'_2$ an *ordered* set by decreeing that $E_1 < E_2$. Note that every member of $\mathcal{E}'_2$ has property $\Pi_{\min}$. Having found the (totally) ordered set $\mathcal{E}'_s \subseteq \mathcal{E}_s$ for some given $s$, we construct the ordered set $\mathcal{E}'_{s+1} \subseteq \mathcal{E}_{s+1}$ in the following iterative manner. Assume that every member of $\mathcal{E}'_s$ has property $\Pi_{\min}$. Let $E \in \mathcal{E}'_s$, and let $F_E$ be the subset of $K$ given by

$$(7.2) \qquad F_E = \Big[\{\langle i,j\rangle : 1 \leq i < j \leq v_{\min}\} + \{\langle v_{\min}, v_{\min}+1\rangle\}\Big] - E.$$

We obtain $\mathcal{E}'_{s+1}$ by adding each such $f \in F_E$ to each such $E$, that is,

$$(7.3) \qquad \mathcal{E}'_{s+1} = \bigcup_{E \in \mathcal{E}'_s} \bigcup_{f \in F_E} \{f + E\}.$$

Note that every member of $\mathcal{E}'_{s+1}$ has property $\Pi_{\min}$. Furthermore, $\mathcal{E}'_{s+1}$ is ordered as follows. For $E', E'' \in \mathcal{E}'_s$ and $f' \in \mathcal{F}_{E'}$, $f'' \in \mathcal{F}_{E''}$, we have $f' + E' < f'' + E''$ if: either $E' < E''$ in $\mathcal{E}'_s$, or $E' = E''$ and $f' < f''$.

It may be seen by an iterative argument involving the quantity $s$ that, for every $E \in \mathcal{E}_s$, there exists $E' \in \mathcal{E}'_s$ such that: there exists a graph isomorphism from $E$ to $E'$ which fixes the edges $e_i$ and $f_i$. It is thus sufficient to consider conjecture instances $(E, e_1, f_1)$, $(E, e_2, f_2)$ for $E$ belonging to some $\mathcal{E}'_s$. Any other conjecture instance is equivalent to some such case. This process of selection reduces the number of conjecture instances, and is effective for small $s$. When $s$ is large, we require the following further economy.

7.2 SKIPPING EQUIVALENT INSTANCES

Having constructed the $\mathcal{E}'_s$, we must verify all conjecture instances $(E, e_i, f_i)$ with $E \in \mathcal{E}'_s$ and $E_i \subseteq E$. However, the cardinalities of the $\mathcal{E}'_s$ grow prohibitively fast despite their selectivity, of the order $\binom{N}{s}$, and in any computer implementation this will tend to exhaust the available resources of memory. In order to address this critical issue, we require a criterion to identify equivalent instances cheaply. This will enable us to further reduce $\mathcal{E}'_s$ by deleting graphs $E$ for which both instances $(E, e_1, f_1)$ and $(E, e_2, f_2)$ are equivalent to earlier ones. We use the following theorem.

**Theorem 7.4.** *Let $E, E' \subseteq K$ and let $e, f \in E$, $e', f' \in E'$ be such that $e \neq f$, $e' \neq f'$. If $j(E, \{e, f\}) = j(E', \{e', f'\})$ then:*

(a) $(E, e, f) \sim (E', e', f')$, and furthermore
(b) for every $H \supseteq E$ there exists $H' \supseteq E'$ such that $(H, e, f) \sim (H', e', f')$.

*Proof.* It suffices to prove part (b), since this contains (a). Let $H \supseteq E$, and recall the definition (7.1) of equivalence. Let $\rho$, $\sigma$, and $\eta$ be as in the proof of Theorem 6.2. We write $X = \{e, f\}$, $X' = \{e', f'\}$. In the present context, $X' = X^\rho$ reads $\{e', f'\} = \{e, f\}^\rho$, so that either $e^\eta = e'$, $f^\eta = f'$, or $e^\eta = f'$, $f^\eta = e'$. We set $H' = H^\eta$. It is now clear that $\eta$ provides one–one correspondences between the pairs $\mathcal{F}^X(H)$, $\mathcal{F}^{X'}(H')$ and $\mathcal{F}_X(H)$, $\mathcal{F}_{X'}(H')$, whence

$$|\mathcal{F}^X(H)| = |\mathcal{F}^{X'}(H')| \quad \text{and} \quad |\mathcal{F}_X(H)| = |\mathcal{F}_{X'}(H')|.$$

Furthermore, if $e^\eta = e'$, $f^\eta = f'$ (respectively, $e^\eta = f'$, $f^\eta = e'$), then $\eta$ provides one–one correspondences between the pairs $\mathcal{F}^e_f(H)$, $\mathcal{F}^{e'}_{f'}(H')$ and $\mathcal{F}^f_e(H)$, $\mathcal{F}^{f'}_{e'}(H')$ (respectively, the pairs $\mathcal{F}^e_f(H)$, $\mathcal{F}^{f'}_{e'}(H')$ and $\mathcal{F}^f_e(H)$, $\mathcal{F}^{e'}_{f'}(H')$). In either case, we have that

$$|\mathcal{F}^e_f(H)| \cdot |\mathcal{F}^f_e(H)| = |\mathcal{F}^{e'}_{f'}(H')| \cdot |\mathcal{F}^{f'}_{e'}(H')|.$$

The condition of (7.1) follows. $\square$

This result is used as follows. Whenever we do the computations required for (5.1) for a given conjecture instance $(E, e, f)$, we place the index $j(E, \{e, f\})$ in a certain database, which we term the *index database*. If later we encounter a conjecture instance $(E', e', f')$ for which $j(E', \{e', f'\})$ lies already in this database, say $j(E', \{e', f'\}) = j(E, \{e, f\})$, we may deduce from Theorem 7.4 that (5.1) is equally valid for $(E', e', f')$ as for $(E, e, f)$. Thus we need not consider the new instance $(E', e', f')$. Indeed, by Theorem 7.4(b), we need consider no instance of the form $(H', e', f')$ for $H' \supseteq E'$.

This observation may be used recursively in order to reduce the ordered sets $\mathcal{E}'_s$ still further. Rather than working with the $\mathcal{E}'_s$ given in (7.3), one works with ordered subsets $\mathcal{E}''_s \subseteq \mathcal{E}'_s$, constructed as follows. We set $\mathcal{E}''_2 = \mathcal{E}'_2$. Suppose that we have found $\mathcal{E}''_2, \mathcal{E}''_3, \ldots, \mathcal{E}''_s$, and we have checked that (5.1) holds for all instances $(E, e_i, f_i)$ with $E \in \mathcal{E}''_t$, $t \leq s$, and all $i$ such that $E_i = \{e_i, f_i\} \subseteq E$. We cycle through the sets $E \in \mathcal{E}''_s$. For each such $E$, we cycle through the edges $f$ belonging to the ordered set $F_E$ of edges defined at (7.2). If both the indices $j(f + E, E_i)$, $i = 1, 2$, lie already in the index database, we do not add the set $f + E$ to $\mathcal{E}''_{s+1}$. If either of the indices do not belong to the database, then we add $f + E$ to $\mathcal{E}''_{s+1}$, and we check any conjecture instance $(f + E, e_i, f_i)$ with $j(f + E, E_i)$ not in the index database. Thus we may write

$$\mathcal{E}''_{s+1} = \bigcup_{E \in \mathcal{E}''_s} \bigcup_{f \in F'_E} \{f + E\},$$

where $F'_E$ is the set of all $f \in F_E$ for which one or both of the above indices do not currently lie in the index database.

**Acknowledgements**

We thank Mark Jerrum, Svante Janson, Malwina Łuczak, Alan Sokal, and Dominic Welsh for discussing this work with us during the 2002 programme entitled 'Computation, Combinatorics and Probability' at the Isaac Newton Institute, Cambridge. The further remarks of Alan Sokal are acknowledged.

Statistical Laboratory, University of Cambridge, Wilberforce Road, Cambridge CB3 0WB, United Kingdom

*E-mail*: g.r.grimmett@statslab.cam.ac.uk, snw22@cam.ac.uk